\documentclass[12pt]{amsart}
\usepackage{latexsym, amsmath, amscd, amssymb, amsthm} 
\usepackage[english]{babel}
 \def\mathbi#1{\textbf{\em #1}}
\begin{document}

\noindent

 \title{ Minimal free resolutions and   Betti diagrams  of the algebras of  $SL_2$-invariants }

\author{Leonid Bedratyuk}

\begin{abstract}
For the algebras of  $SL_2$-invariants of small homological dimension  theirs  minimal free graded resolutions and   graded Betti diagrams calculated.
\end{abstract}
\maketitle

{\bf 1.} Let $K$ be a field, ${\rm char}K=0.$ Let $V_d $ be  $d+1$-dimensional  $SL_2$-module of binary forms of degree $d$ and let $V_{\mathbi{d}} =
V_{d_1} \oplus V_{d_2} \oplus \cdots \oplus V_{d_n},$ ${\mathbi{d}=(d_1,d_2,\ldots, d_n).}$ Denote by  $K[V_{\mathbi{d}}]^{SL_2}$   the algebra of polynomial $SL_2$-invariant functions on $V_{\mathbi{d}}.$  It is well known that the algebra
 $\mathcal{I}:=K[V_{\mathbi{d}}]^{SL_2}$ is finitely generated  and let $f_1,f_2,\ldots,f_m$  be  its minimal generating set.  The measure of the
intricacy of this algebra is the length of their chains of syzygies, called
homological(or  proective)  dimension ${\rm hd} \, \mathcal{I}.$  In \cite{Pop} Popov gave  a classification of the cases in which ${\rm hd} \, \mathcal{I} \leq 10 $ for a single binary form ($n = 1$) or ${\rm hd}\,  \mathcal{I} \leq 3 $ for $n>1.$ Recently Brouwer and Popoviciu \cite{Bro} extend the results and determine for $n = 1$ the cases with
${\rm hd} \, \mathcal{I} \leq 100,$ and for $n > 1$ those with ${\rm hd}\,  \mathcal{I} \leq 15.$

 In this short  notes we present results of calculations  of    the free minimal graded resolutions and the graded Betti diagrams for  the algebras of $SL_2$-invariants $\mathcal{I}$ in the cases  ${\rm hd}\, \mathcal{I} \leq 8.$    

{\bf 2.}
Let $R := K[x_1,\ldots , x_m]$ be positively graded by $\deg(x_i)=\deg(f_i),$ $ i=1,\ldots,m.$  Recall that a finite graded free resolution of $\mathcal{I}$  of length $l={\rm hd}\, \mathcal{I} $ is an exact sequence of  $R$-modules
 $$
0 \longrightarrow  F_{l} \longrightarrow  F_{l-1} \longrightarrow \cdots \longrightarrow F_{1} \longrightarrow F_0 \longrightarrow \mathcal{I} \longrightarrow 0,
$$
where $F_i=\oplus_{j}R(-j)^{\beta_{i,j}},$ $F_0=R$  are finitely generated graded free $R$-module.
The image $F_i$  is called the $i$-th module of syzygies of $\mathcal{I}.$  Since  the algebra $\mathcal{I}$ is Cohen-Macaulay then  the   Auslander-Buchsbaum  theorem implies that $l=\dim \mathcal{I}  - {\rm tr}\deg_K \mathcal{I},$  see also \cite{Pop}.

  The numbers $\beta_{i,j}$  are  called the graded Betti numbers which  can be arranged into a graded  Betti diagram $\beta(\mathcal{I})$:

\begin{center}
 \begin{tabular}{c|cccccc}
  $-j\backslash i$  &0 &  1           & 2&  3 & \ldots & $l$  \\ 
\hline 
   0      &1& -&-&- & \ldots& \\
    1      &-& $ \beta_{1,1} $& $\beta_{2,1}$& $ \beta_{3,1}$ & \ldots& \\
    
    2     &-& $\beta_{1,2}$ & $\beta_{2,2}$&  $\beta_{3,2}$ & \ldots & \\
    3     &-& $\beta_{1,3}$ & $\beta_{2,3}$&  $\beta_{3,3}$ & \ldots & \\
    \vdots     &  \vdots &  \vdots &   \vdots &  \vdots & \\
    $j^*$      &-& -&-&- & \ldots&1\\
\hline 

\end{tabular}
\end{center}
\vspace{0.3cm}
where $j^*$  is maximal of the indexes $j$ appears in  the resolutions.

The Hilbert-Poincar\'e series of the algebra $\mathcal{I}$  can be  recovered from the Betti  numbers by
$$
\mathcal{P}(\mathcal{I},z)=\frac{\displaystyle \sum_{i=0}^l\sum_{j \in \mathbb{Z}} (-1)^j \beta_{i,j} z^j}{\displaystyle \prod_{i=1}^{m}(1-z^{\deg(f_i)})}
$$
The  inverse statement is false, see below the   case $\mathcal{I}=k[2V_1 \oplus V_3]^{SL_2}.$

To calculate the minimal free  resolutions we used the following algorithm:

\begin{itemize}
	\item find a minimal generating set  of the algebra of invariants $\mathcal{I}$  by using our Maple-package {\tt SL\_2 invariants}, see \cite{Bed};
	
	\item find a minimal generating set of the first syzygy $R$-module $F_1$ by using {\rm CoCoA,}  see \cite{CA};
	
	\item find a minimal free  resolution  of $R/F_1$ by using {\rm CoCoA.}
\end{itemize}

For  example, let us calculate the free minimal  resolution of $\mathcal{I}:=k[ 3 V_1 \otimes V_2]^{SL_2}.$
By using the  Maple  package procedure ${\tt Min\_Gen\_Set\_Invariants([1,1,1,2])}$ we get the  following 10 elements of the  minimal generating set  of $\mathcal{I}:$
 \begin{gather*} L=[ -x_{{0}}u_{{1}}v_{{1}}+x_{{0}}u_{{0}}v_{{2}}-x_{{1}}u_{{0}}v_{{1}}+x_{{1}}u_{{1}}v_{{0}},-2\,u_{{0}}u_{{1}}v_{{1}}+{u_{{0}}}^{2}v_{
{2}}+{u_{{1}}}^{2}v_{{0}},-2\,{v_{{1}}}^{2}+2\,v_{{0}}v_{{2}},\\ -{x_{{1}
}}^{2}v_{{0}}-{x_{{0}}}^{2}v_{{2}}+2\,x_{{0}}x_{{1}}v_{{1}},-x_{{0}}u_
{{1}}+x_{{1}}u_{{0}},y_{{1}}u_{{1}}v_{{0}}-y_{{1}}u_{{0}}v_{{1}}-y_{{0
}}u_{{1}}v_{{1}}+y_{{0}}u_{{0}}v_{{2}},\\{y_{{1}}}^{2}v_{{0}}-2\,y_{{0}}
y_{{1}}v_{{1}}+{y_{{0}}}^{2}v_{{2}},-x_{{0}}y_{{1}}+x_{{1}}y_{{0}},-y_
{{0}}u_{{1}}+y_{{1}}u_{{0}},x_{{0}}y_{{0}}v_{{2}}-x_{{0}}y_{{1}}v_{{1}
}-x_{{1}}y_{{0}}v_{{1}}+x_{{1}}y_{{1}}v_{{0}}].
\end{gather*}
By  using CocoA  commands\\
 
\noindent
{\tt Use QQ[x[0..1],y[0..1],u[0..1],v[0..2]];\\ 
SAM := SubalgebraMap(L);K:=Ker(SAM);K;} \\

\noindent
  we  calculate the kernel $J$ of map 
$\varphi: k[x_1,\ldots, x_{10}] \to \mathcal{I}, x_i \mapsto L[i]:$

{\tt J:= Ideal($- x_{5}x_{7} + x_{6}x_{8} + x_{9}x_{10},
  x_{5}x_{6} - x_{2}x_{8} - x_{1}x_{9},
   - x_{1}x_{8} + x_{4}x_{9} + x_{5}x_{10},
   - x_{6}^{2} + x_{2}x_{7} -1/2x_{3}x_{9}^{2},
   - x_{4}x_{6} -1/2x_{3}x_{5}x_{8} - x_{1}x_{10},
   - x_{1}x_{6} -1/2x_{3}x_{5}x_{9} + x_{2}x_{10},
   - x_{1}x_{7} -1/2x_{3}x_{8}x_{9} + x_{6}x_{10},
   - x_{1}^{2} - x_{2}x_{4} -1/2x_{3}x_{5}^{2},
   - x_{4}x_{7} -1/2x_{3}x_{8}^{2} - x_{10}^{2},
  x_{2}x_{4}x_{8} + 1/2x_{3}x_{5}^{2}x_{8} + x_{1}x_{4}x_{9} + x_{1}x_{5}x_{10})$} \\

Then by  using  CoCoA commands \\

\noindent
{\tt
Use R::=QQ[x[1..10]],Weights(3, 3, 2, 3, 2, 3, 3, 2, 2, 3);\\
Res(R/J);
}\\

\noindent
we  get the  minimal free  resolution of the ring $R/J:$

$0
 \longrightarrow R(-17)
 \longrightarrow R(-11)^{6} \oplus R(-12)^{3}
 \longrightarrow R(-8)^{8} \oplus R(-9)^{8}
 \longrightarrow R(-5)^{3} \oplus R(-6)^{6}
 \longrightarrow R$

 Since $R/J \cong \mathcal{I}$  we  may extend it to the  minimal free  resolution of $\mathcal{I}:$
 
   $0
 \longrightarrow R(-17)
 \longrightarrow R(-11)^{6} \oplus R(-12)^{3}
 \longrightarrow R(-8)^{8} \oplus R(-9)^{8}
 \longrightarrow R(-5)^{3} \oplus R(-6)^{6}
 \longrightarrow R \longrightarrow \mathcal{I} \longrightarrow 0.$
 
 The graded Betti diagram for $\mathcal{I}$ is as following

\begin{center}

  \begin{tabular}{c|ccccc}
  $-j\backslash i$ &0  &  1  & 2&  3 & 4   \\ 
\hline 
    0  &1 & - & -&  - & - \\
    5  &- & 3 & -&  - & - \\
    
    6   &-  & 6 & -&  - & - \\
    8   &-  & - & 8&  - & - \\
    9    &- & - & 8&  - & - \\
    11  &-  & - &- &  6 & -  \\
    12  &-   & - & -& 3 & -  \\
    17   &-  & - & -&  - & 1\\
\hline 

\end{tabular}
\end{center}

\vspace{0.5cm}

{\bf 3.}  Below we list  the minimal free graded resolution and the graded Betti diagrams of the algebras  of invariants  in the cases ${\rm hd}\, \mathcal{I} \leq 8.$

\subsection{{\rm hd}\, $\mathcal{I}=0$} In this  case we have  $ \mathcal{I} \in \{V_1^{SL_2}
, V_2^{SL_2}, V_3^{SL_2}, V_4^{SL_2}, (2V_1)^{SL_2}, (V_1 \oplus V_2)^{SL_2},\\ (2V_2)^{SL_2}, (3V_1)^{SL_2} \}$ and the   minimal free  resolutions have the following form:
$$0
  \longrightarrow R \longrightarrow \mathcal{I} \longrightarrow 0;$$

\subsection{{\rm hd}\, $\mathcal{I}=1$}  

$\mathcal{I}:=k[V_5]^{SL_2},$ 

 $$0
 \longrightarrow R(-36)
 \longrightarrow R \longrightarrow \mathcal{I} \longrightarrow 0.$$

$\mathcal{I}:=k[V_6]^{SL_2},$  $$0
 \longrightarrow R(-30)
 \longrightarrow R \longrightarrow \mathcal{I} \longrightarrow 0;$$
 
 The  above cases $d=5,6$ are well known classical results, see \cite{Hilb}.

 $\mathcal{I}:=k[V_1 \oplus  V_3]^{SL_2},$  $$0
 \longrightarrow R(-12)
 \longrightarrow R  \longrightarrow \mathcal{I} \longrightarrow 0 ;$$

 $\mathcal{I}:=k[V_1 \oplus  V_4]^{SL_2},$  $$0
 \longrightarrow R(-18)
 \longrightarrow R  \longrightarrow \mathcal{I} \longrightarrow 0  ;$$

 $\mathcal{I}:=k[ V_2 \oplus  V_3]^{SL_2},$  $$0
 \longrightarrow R(-14)
 \longrightarrow R  \longrightarrow \mathcal{I} \longrightarrow 0 ;$$
 
   $\mathcal{I}:=k[ V_2 \oplus  V_4]^{SL_2},$  $$0
 \longrightarrow R(-12)
 \longrightarrow R  \longrightarrow \mathcal{I} \longrightarrow 0 ;$$
 
  $\mathcal{I}:=k[ V_4 \oplus  V_4]^{SL_2},$  $$0
 \longrightarrow R(-12)
 \longrightarrow R  \longrightarrow \mathcal{I} \longrightarrow 0 ;$$
 
  $\mathcal{I}:=k[V_1 \oplus   V_1 \oplus  V_2]^{SL_2},$  $$0
 \longrightarrow R(-6)
 \longrightarrow R  \longrightarrow \mathcal{I} \longrightarrow 0 ;$$

 $\mathcal{I}:=k[V_1 \oplus   V_2 \oplus  V_2]^{SL_2},$  $$0
 \longrightarrow R(-8)
 \longrightarrow R  \longrightarrow \mathcal{I} \longrightarrow 0 ;$$

 $\mathcal{I}:=k[V_2 \oplus   V_2 \oplus  V_2]^{SL_2},$  $$0
 \longrightarrow R(-6)
 \longrightarrow R  \longrightarrow \mathcal{I} \longrightarrow 0 ;$$
 
 $\mathcal{I}:=k[V_1 \oplus   V_1 \oplus  V_1 \oplus  V_1]^{SL_2},$  $$0
 \longrightarrow R(-4)
 \longrightarrow R  \longrightarrow \mathcal{I} \longrightarrow 0 ;$$
 
 \subsection{{\rm hd}\, $\mathcal{I}=2$}

 $\mathcal{I}:=k[V_3  \oplus  V_3]^{SL_2},$\\  $0  \longrightarrow R(-20) 
 \longrightarrow R(-8) \oplus R(-12)
 \longrightarrow R  \longrightarrow \mathcal{I} \longrightarrow 0 ;$

\begin{center}
 \begin{tabular}{c|ccc}
  $-j\backslash i$ &0  &  1  & 2    \\ 
\hline 
  0 &1   & - & -  \\
    8 &-  & 1 & -  \\
    
    12 &- & 1 & -   \\
    20  &-  & - & 1  \\
    
\hline 

\end{tabular}
\end{center}
 
  \subsection{{\rm hd}\, $\mathcal{I}=3$}
 
 \vspace{0.5cm}
 $\mathcal{I}:=k[V_8]^{SL_2},$
 
   $0  \longrightarrow R(-45) \longrightarrow R(-25) \oplus R(-26)\oplus R(-27)\oplus R(-28)\oplus R(-29)
 \longrightarrow R(-16) \oplus R(-17)\oplus R(-18)\oplus R(-19)\oplus R(-20)
  \longrightarrow R \longrightarrow \mathcal{I} \longrightarrow 0 ;$
  
  \begin{center}
 \begin{tabular}{c|cccc}
  $-j\backslash i$ &0&  1  & 2 & 3   \\ 
\hline 
     0  &1 & -   & - &-  \\
    16  &- & 1   & - &-  \\
    
    17  &- & 1   & - &-  \\
    18   &-& 1 & - &- \\
    19   &-& 1 & - &- \\
    
    20  &-& 1 & - & - \\
    25   &- & - & 1 &- \\
    26  &- & - & 1 & -\\
    
    27 &- & - & 1 & - \\
    28  &- & - & 1 &- \\
    29  &- & - & 1 &- \\
    
   45 &- & - & - &1  \\
    
\hline 

\end{tabular}
\end{center}
  
  It is result of Shioda \cite{Shi}.\\

   $\mathcal{I}:=k[5 V_1]^{SL_2},$
   
   $0  \longrightarrow R(-10) \longrightarrow R(-6)^5 
 \longrightarrow  R(-4)^5 
  \longrightarrow R \longrightarrow \mathcal{I} \longrightarrow 0 ;$
  
    \begin{center}
 \begin{tabular}{c|cccc}
  $-j\backslash i$ &0 &  1  & 2 & 3   \\ 
\hline
    0   &1& -  & - &-  \\
    4   &- & 1   & - &-  \\    
    6   &-& -   & 1 &-  \\
    10  &-& - & - &1 \\
    
\hline 

\end{tabular}
\end{center}
 
  \subsection{{\rm hd}\, $\mathcal{I}=4$}   $\mathcal{I}:=k[3V_1 \oplus V_2]^{SL_2},$
 
$0 \longrightarrow R(-17) \longrightarrow R(-11)^6 \oplus R(-12)^3 \longrightarrow R(-8)^8 \oplus R(-9)^8 \longrightarrow R(-5)^3 \oplus R(-6)^6 \longrightarrow  R \longrightarrow \mathcal{I} \longrightarrow 0$

  \begin{center}
 \begin{tabular}{c|ccccc}
  $-j\backslash i$ &0&  1  & 2 & 3&4   \\ 
\hline
    0  &1& -   & - &-&-  \\
    5   &-& 3   & - &-&-  \\
    
    6  &-& 6   & - &- &- \\
    8   &-& - & 8 &-&-\\
    9   &-& - & 8 &- &- \\
    
    11  &-& - & - & 6 &- \\
    12   &- & - & - &3&-\\
    17  &- & - & - &-& 1\\
            
\hline 

\end{tabular}
\end{center}
\subsection{{\rm hd}\, $\mathcal{I}=5$}   $\mathcal{I}:=k[V_1 \oplus 3V_2]^{SL_2},$

$0 \longrightarrow R(-25) \longrightarrow R(-17)^6 \oplus R(-18)^4 \oplus R(-19)^4 \longrightarrow R(-13)^8 \oplus R(-14)^{12} \oplus R(-15)^{12} \oplus R(-16)^3 \longrightarrow R(-9)^3 \oplus R(-10)^{12} \oplus R(-11)^{12} \oplus R(-12)^8 \longrightarrow R(-6)^4 \oplus R(-7)^4 \oplus R(-8)^6 \longrightarrow R \longrightarrow \mathcal{I} \longrightarrow 0 $

\begin{center}
 \begin{tabular}{c|cccccc}
  $-j\backslash i$ &0&  1  & 2 & 3&4&5   \\ 
\hline 
    0   &1& -   & - &-&-&-  \\
    6   &-& 4   & - &-&-&-  \\
    7   &-& 4   & - &-&-&-  \\    
    8   &-& 6   & - &-&-&-  \\
    9   &-& -   & 3 &-&-&-  \\
    10  &-& -   & 12 &-&-&-  \\    
    11   &-& -   & 12 &-&-&-  \\
    12   &-& -   & 8 &-&-& - \\
    13   &-& -   & - &8&-&-  \\    
    14   &-& -   & - &12&-&-  \\ 
    15   &-& -   & - &12&-&-  \\
    16   &-& -   & - &3&-&-  \\    
    17   &-& -   & - &-&6&-  \\ 
    18   &-& -   & - &-&4&-  \\
    19   &-& -   & - &-&4&-  \\    
    25   &-& -   & - &-&-&1  \\ 
\hline 

\end{tabular}
\end{center}

$\mathcal{I}:=k[4V_2]^{SL_2},$

$0 \longrightarrow R(-20) \longrightarrow R(-14)^{10} \oplus R(-15)^4 \longrightarrow R(-11)^{20} \oplus R(-12)^{15} \longrightarrow R(-8)^{15} \oplus R(-9)^{20} \longrightarrow  R(-5)^4 \oplus R(-6)^{10} \longrightarrow R \longrightarrow \mathcal{I} \longrightarrow 0 $

\begin{center}
 \begin{tabular}{c|cccccc}
  $-j\backslash i$ &0&  1  & 2 & 3&4&5   \\ 
\hline 
    0   &1& -   & - &-&-&-  \\
    5   &-& 4   & - &-&-&-  \\
    6  &-& 10   & - &-&-&-  \\    
    8   &-& -   & 15 &-&-&-  \\
    9   &-& -   & 20 &-&-&-  \\
    11  &-& -   & - &20&-&-  \\    
    12   &-& -   & - &15&-&-  \\
    14   &-& -   & - &-&10& - \\
    15   &-& -   & - &-&4&-  \\    
    20   &-& -   & - &-&-&1  \\ 
    \hline 

\end{tabular}
\end{center}

\subsection{{\rm hd}\, $\mathcal{I}=6$}   $\mathcal{I}:=k[2V_1 \oplus  2 V_2]^{SL_2},$

$0
 \longrightarrow R(-28)
 \longrightarrow R(-20)^{6} \oplus R(-21)^{8} \oplus R(-22)^{6}
 \longrightarrow R(-16)^{8} \oplus R(-17)^{24} \oplus R(-18)^{24} \oplus R(-19)^{8}
 \longrightarrow R(-12)^{3} \oplus R(-13)^{24} \oplus R(-14)^{36} \oplus R(-15)^{24} \oplus R(-16)^{3}
 \longrightarrow R(-9)^{8} \oplus R(-10)^{24} \oplus R(-11)^{24} \oplus R(-12)^{8}
 \longrightarrow R(-6)^{6} \oplus R(-7)^{8} \oplus R(-8)^{6}
 \longrightarrow R \longrightarrow \mathcal{I} \longrightarrow 0$
 
 \begin{center}
 \begin{tabular}{c|ccccccc}
  $-j\backslash i$ &0&  1  & 2 & 3&4&5&6   \\ 
\hline 
    0   &1& -   & - &-&-&-&-  \\
    5   &-& 4   & - &-&-&-&-  \\
    6  &-& 10   & - &-&-&-&-  \\    
    8   &-& -   & 15 &-&-&-&-  \\
    9   &-& -   & 20 &-&-&-&-  \\
    11  &-& -   & - &20&-&-&-  \\    
    12   &-& -   & - &15&-&-&-  \\
    14   &-& -   & - &-&10& -&- \\
    15   &-& -   & - &-&4&-&-  \\    
    20   &-& -   & - &-&-&1&-  \\ 
    \hline 

\end{tabular}
\end{center}

 $\mathcal{I}:=k[6V_1 ]^{SL_2},$
  
 $0
 \longrightarrow R(-18)
 \longrightarrow R(-14)^{15}
 \longrightarrow R(-12)^{35}
 \longrightarrow R(-8)^{21} \oplus R(-10)^{21}
 \longrightarrow R(-6)^{35}
 \longrightarrow R(-4)^{15}
 \longrightarrow R \longrightarrow \mathcal{I} \longrightarrow 0 $
 
 \begin{center}
 \begin{tabular}{c|ccccccc}
  $-j\backslash i$ &0&  1  & 2 & 3&4&5&6   \\ 
\hline 
    0   &1& -   & - &-&-&-&-  \\
    4   &-& 15   & - &-&-&-&-  \\
    6   &-& -    & 35 &-&-&-&-  \\    
    8   &-& -    &   &21&-&-&-  \\
    10  &-& -   & -  &21&-&-&-  \\
    12  &-& -   & -  &-&35&-&-  \\    
    14   &-& -   & - &-&-&15&-  \\
    18   &-& -   & - &-&-& -&1 \\
        \hline 

\end{tabular}
\end{center}
 \subsection{{\rm hd}\, $\mathcal{I}=8$}   $\mathcal{I}:=k[2V_1 \oplus V_3]^{SL_2},$
 
 $0 \longrightarrow R(-50) 
 \longrightarrow R(-38)^{10} \oplus R(-40)^{15} \oplus R(-42)^{10} 
 \longrightarrow R(-32)^{20} \oplus R(-34)^{60} \oplus R(-36)^{60} \oplus R(-38)^{20} 
 \longrightarrow R(-26)^{15} \oplus R(-28)^{90} \oplus R(-30)^{140} \oplus R(-32)^{90} \oplus R(-34)^{15} 
  \longrightarrow R(-20)^4 \oplus R(-22)^{60} \oplus R(-24)^{160} \oplus R(-26)^{160} \oplus R(-28)^{60} \oplus R(-30)^4
   \longrightarrow R(-16)^{15} \oplus R(-18)^{90} \oplus R(-20)^{140} \oplus R(-22)^{90} \oplus R(-24)^{15} 
  \longrightarrow R(-12)^{20} \oplus R(-14)^{60} \oplus R(-16)^{60} \oplus R(-18)^{20} 
  \longrightarrow R(-8)^{10} \oplus R(-10)^{15} \oplus R(-12)^{10}
  \longrightarrow R \longrightarrow \mathcal{I} \longrightarrow 0$
 
 \begin{center}
 \begin{tabular}{c|ccccccccc}
  $-j\backslash i$ &0&  1  & 2 & 3&4&5&6&7&8   \\ 
\hline 
    0   &1& -   & - &-&-&-&-&-&-  \\
    8   &-& 10   & - &-&-&-&-&-&-  \\
    10   &-& 15  & - &-&-&-&-&-&-  \\
    12   &-& 10  & 20 &-&-&-&-&-&-  \\
    14   &-& -   & 60 &-&-&-&-&-&-  \\
    16   &-& -   & 60 &15&-&-&-&-&-  \\
    18   &-& -   & 20 &90&-&-&-&-&-  \\
    20   &-& -   & -  &140&4&-&-&-&-  \\
    22   &-& -    & - &90&60&-&-&-&-  \\
    24   &-& -   & -  &15&160&-&-&-&-  \\
    26   &-& -   & -  &- &160&15&-&-&-  \\
    28   &-& -   & -  &- &60 &90&-&-&-  \\
    30   &-& -   & -  &- &4  &140&-&-&-  \\
    32  &-& -   & -   &- &-  &90&20&-&-  \\
    34  &-& -   & -   &- &-  &15&60&-&-  \\
    36   &-& -   & -  &- &-  &- &60&-&-  \\
    38   &-& -   & -  &- &-  &- &20&10&-  \\
    40   &-& -   & -  &- &-  &- &- &15&-  \\
    42  &-& -   & -   &- &-  & -&- &10&-  \\
    50  &-& -   & -   &- &-  &- &- &-&1  \\
    \hline 
  
\end{tabular}
\end{center}

{\bf Conjecture.} For algebras of $SL_2$-invariants $\mathcal{I}$ its graded Betti diagram is palindromic, i.e.  $\beta_{l-i,j^*-j}=\beta_{i,j}.$

\end{document}